\documentclass[letterpaper, 10 pt, twocolumn, swedish]{ieeetran}

\usepackage{graphicx}
\usepackage[english]{babel}
\usepackage{amsmath}
\usepackage{amsfonts}
\usepackage{dansimath}
\usepackage{cite}
\usepackage{stfloats}

\title{Stability analysis of Model Predictive Controllers \\ using Mixed Integer Linear Programming}

\author{Daniel Simon and Johan L\"ofberg \thanks{The authors are with the Departement of Electrical Engineering, Link\"oping University, SE-581 83 Link\"oping, Sweden, e-mail: dansi@isy.liu.se, johanl@isy.liu.se}
\thanks{This work has been submitted to the 55\textsuperscript{th}, IEEE Conference on Decision
  and Control 2016}}

\begin{document}
\maketitle

\begin{abstract}
It is a well known fact that finite time optimal controllers, such as MPC does not necessarily result in closed loop stable systems. Within the MPC community it is common practice to add a final state constraint and/or a final state penalty in order to obtain guaranteed stability. However, for more advanced controller structures it can be difficult to show stability using these techniques. Additionally in some cases the final state constraint set consists of so many inequalities that the complexity of the MPC problem is too big for use in certain fast and time critical applications. In this paper we instead focus on deriving a tool for a-postiori analysis of the closed loop stability for linear systems controlled with MPC controllers. We formulate an optimisation problem that gives a sufficient condition for stability of the closed loop system and we show that the problem can be written as a \emph{Mixed Integer Linear Programming Problem} (MILP).
\end{abstract}
 
\section{Introduction}
A linear Model Predictive Controller (MPC) solves, online in each sample instant, a finite time horizon optimal control problem of the form
\begin{subequations}  
\begin{align}
V_k^* = \minim{x_{k+i},u_{k+i}} \;\; & \sum_{i=0}^{N-1} \ell(x_{k+i},u_{k+i}) + \Psi(x_{k+N})\label{eq:MPC_obj}\\
\st \;\;\; & x_{k+i+1} = Ax_{k+i} + Bu_{k+i} \label{eq:MPC_dyn}\\
 & Ex_{k+i} \leq f \quad i = \range{1}{N-1} \label{eq:MPC_state}\\
 & Tx_{k+N} \leq t \label{eq:MPC_final}\\
 & Gu_{k+i} \leq h \quad i = \range{1}{N-1} \label{eq:MPC_control}
\end{align} \label{eq:MPC}
\end{subequations} 
and implements the optimal solution, $u_k^*$, as the input to the system in a receding horizon fashion. 

While the system dynamics \eqref{eq:MPC_dyn}, the state and control constraints \eqref{eq:MPC_state} and \eqref{eq:MPC_control} reflects the system and the requirements on it, the final state constraint \eqref{eq:MPC_final} and penalty, $\Psi(x_{k+N})$ are in a sense artificial and added merely to ensure recursive feasibility and stability. 

Stability of the MPC control law \eqref{eq:MPC} is most often proven a-priori by showing that the objective function \eqref{eq:MPC_obj} is a valid Lyapunov function for the closed loop system by designing $\Psi(x_{k+N})$ and the terminal set such that 
\begin{equation}
V^*_{k+1} - V^*_{k} 
\label{eq:dV}
\end{equation}
is guaranteed to be less that zero. See, e.g., \cite{Mayne2000} for details.

However there exist many MPC formulations, such as move blocking \cite{Cagienard2007} or soft constraints \cite{Zeilinger2014} (i.e., slack), for which it can be difficult to show stability using this standard framework presented in \cite{Mayne2000}.  Even if it allows for the possibility to guarantee stability, the addition of the terminal constraints can add a significant complexity to the original problem which might be unnecessary and limit the applicability of MPC within certain fields. In reference tracking MPC the terminal set can become very complex and some times not even finitely determined, see e.g., \cite{Chisci2003} or \cite{Limon2008}. Therefore it can often be beneficial to analyse and verify the stability of a certain design rather than building in the stability by adding extra constraints.

The problem of analysing stability of optimisation based controllers is no new field. In \cite{Primbs2000} the authors derive a stability test based on bounds of the cost function.  The author of \cite{Primbs2001} uses the KKT conditions for the MPC controller \eqref{eq:MPC} to derive, using the S-procedure, an LMI which gives a sufficient condition for stability. Several papers have followed up on this idea such as, e.g., \cite{Ahmad2014,Li2005,Li2006,Lovaas2007,Korda2015}. These papers have extended the idea to hold for more general cases or improved the complexity of the resulting LMI. E.g., in \cite{Korda2015} the authors extends the ideas in \cite{Primbs2001} to hold for more general systems and show that given that the system and its constraints are polynomial, the stability can be analysed using sum-of-squares programming. In \cite{Li2005} the authors derive an LMI of much lower dimension than that of \cite{Primbs2001} which improve on the complexity of the problem. 

In this paper we also exploit the KKT conditions of the MPC problem \eqref{eq:MPC} to analyse the difference of the value function \eqref{eq:dV}. But instead of deriving an LMI condition we formulate the problem as an indefinite quadratic bilevel optimisation problem. We then show that this problem can be rewritten as a Mixed Integer Linear Program (MILP), which can have major computational advantages compared to using the LMI formulation.
 
\section{The MILP Stability test} \label{seq:the_algorithm}
In this section we will derive the proposed stability test. We will formulate an optimisation problem that uses the MPC problem's objective function as a candidate Lyapunov function and then minimises the difference of the Lyapunov function between two consecutive time steps. Since the optimisation problem only test the validity of a certain Lyapunov function candidate the formulated problem can only verify stability, not prove instability. Therefore it is a sufficient but not necessary condition for stability. Nevertheless, in situations where the value function is non-decreasing at some point, it might be an indication of problems in the design and furtrher investigations and simulations needs to be made.

The resulting optimisation problem is an indefinite quadratic bilevel optimisation problem which are very difficult to solve. We will finally show how this indefinite problem can be rewritten as a mixed integer linear programming problem. 

If the objective function \eqref{eq:MPC_obj} is a quadratic or LP-representable function in $x$ and $u$ we can rewrite the MPC problem \eqref{eq:MPC}, at time step $k$, in a more compact form as
\begin{align}
\minimize{U_{k}} \;\; & V_k\label{eq:MPC2}\\
\st  \;\; & Ex_{k} + FU_{k} \leq b \nonumber
\end{align}
where \[V_k = \half U_{k}^THU_{k} + U_{k}^TGx_{k} + \half x_{k}^T\bar{Q}x_{k}\] and where $U_k = [u^T_k, u^T_{k+1}, \dots, u^T_{k+N-1}]^T$, $x_k$ is the current measured state and the matrices, $H$, $G$, $\bar{Q}$, $E$, $F$ and $b$ are suitably defined. The objective function, $V_{k+1}$, at time $k+1$ is analogously defined. 

Let us now formulate the stability test as the following optimisation problem 
\begin{subequations} \label{eq:stabilitytest}
\begin{align} 
  \minimize{U^*_k,U^*_{k+1},x_k} & V^*_k - V^*_{k+1} \label{eq:stab_obj}\\
 \st \;\; & x_{k+1} = Ax_k + Bu^*_k \\ 
& U_k^* = \argmin{} V_k \label{eq:mpc_in_stab} \\
& U_{k+1}^* = \argmin{} V_{k+1} \label{eq:mpc_in_stab_2} 
\end{align}
\end{subequations}
Hence we want to find the state $x_k$ of the system that result in the smallest possible difference in our candidate Lyapunov function when controlled with the MPC controller.

If this difference is positive (or equal to zero) this means that $V_k$ is a valid Lyapunov function for the system and hence it is stable. On the other hand if this difference is negative this means that we have an increase in the Lyapunov function candidate for som point $x_k$ and hence it is not a valid Lyapunov function for the system. 

Note that it is straightforward to extend the algorithm to incorporate a condition of sufficient decrease and thus asymptotic stability and robustness, e.g., by modifying the objective function to  
\[ V^*_{k} - V^*_{k+1} - \epsilon \left(x_k^TQx_k - u_k^TRu_k\right) \geq 0\]

At this point it should be pointed out that it is assumed, throughout the paper, that the MPC algorithms are recursively feasible. This must of course be tested before one can apply the proposed stability test. A method for performing the feasibility test is presented in \cite{Lofberg2012}.

Note also that in the case of the MPC problem having several non-unique solutions then the stability test derived here is a pessimistic bound on the stability, i.e., it selects the worst case combinations of optimal points. 

Since the MPC problem \eqref{eq:MPC2} is a convex QP we can replace it in \eqref{eq:mpc_in_stab} and \eqref{eq:mpc_in_stab_2} with the necessary and sufficient KKT conditions 
\begin{subequations} \label{eq:KKT}
\begin{align} %
  HU_k + Gx_k + F^T\lambda_k &= 0 \\
 Ex_k + FU_k -b &\leq 0 \\
\lambda_k &\geq 0 \\
\lambda^T(Ex_k + FU_k -b) &= 0 \label{eq:complement}
\end{align}
\end{subequations}  
and similar for time $k+1$.

Note that equations \eqref{eq:complement} are bilinear constraints but they can for each row be modeled using a \emph{Big-M} reformulation as four linear constraints with a binary variable, $z^{(i)}_k$, as
\begin{multline} \lambda^{(i)}_k(e_i^Tx_k + f_i^TU_k - b_i) = 0 \Rightarrow \\
0 \leq \lambda^{(i)}_k \leq m^{(i)}_1z^{(i)}_k, \\ 0 \leq b_i-e_i^Tx_k - f_i^TU_k \leq m^{(i)}_2(\ones-z^{(i)}_k) \label{eq:compl_bin}
\end{multline}
We can see that the binary variable $z_k^{(i)}$ forces either the constraint, $b_i-e_i^Tx_k - f_i^TU_k$, to be equal to zero or the dual variable, $\lambda_k^{(i)}$, to be equal to zero, ensuring that the complementarity constraint holds. In other words, we can view this as that $z_k^{(i)}$ encodes whether the constraint is active or not, i.e., $z_k^{(i)}=1 \Rightarrow e_i^Tx_k + f_i^TU_k = b_i$.

Using the KKT conditions \eqref{eq:KKT} and the binary reformulation \eqref{eq:compl_bin} we can write the problem \eqref{eq:stabilitytest} using the notation, ${y = [U^T_k,x^T_k,\lambda^T_k,U^T_{k+1},x^T_{k+1},\lambda^T_{k+1}]^T}$ and  $\bar{z}=[z^T_k,z^T_{k+1}]^T$ as
\begin{subequations} \label{eq:blocked}
\begin{align}
  \minimize{y,\bar{z}} & \half y^T\bar{H}y \\
  \st \;\; & \bar{E}y = 0 \\
  & \bar{A}y \leq \bar{b} + \bar{d}\bar{z} \label{eq:blocked_constr}
\end{align}  
\end{subequations}

where
\[ \bar{H} = \begin{bmatrix} H & G & 0 & 0&0&0\\ G^T & \bar{Q} & 0 & 0&0&0 \\ 
0&0&0&0&0&0 \\
0 & 0 & 0 & -H & -G & 0  \\ 
0 & 0 & 0 & -G^T & -\bar{Q} & 0  \\
0&0&0&0&0&0 \end{bmatrix} \]

\[ \bar{E} = \begin{bmatrix} \bar{B} & A & 0 & 0 & -I & 0  \\
H & G & F^T & 0 & \dots \\
 & \dots & 0 & H & G & F^T  \end{bmatrix} \]

\[ \bar{A} = \begin{bmatrix} F & E & 0 & \dots \\
-F & -E & 0 & 0 & \dots \\
\dots & 0 & -I &  0 & \dots \\
\dots & 0 & I &  0 & \dots \\
&\dots&0& F & E & 0  \\
&\dots&0& -F & -E & 0 \\
&\dots&0& 0 & 0 & -I \\
&\dots&0& 0 & 0 & I \end{bmatrix} \]

\[ \bar{b} = \begin{bmatrix} b \\ M_1\ones - b \\ 0 \\ 0 \\ b \\ M_1\ones - b \\ 0 \\ 0 \end{bmatrix}, \quad  \bar{d} = \begin{bmatrix} 0 & 0 \\ -M_1 & 0 \\ 0 & 0 \\ M_2 & 0 \\0 & 0 \\  0 & -M_1  \\ 0 & 0 \\ 0 & M_2 \end{bmatrix} \]

Even though we have eliminated the bilinear constraints the problem is still an indefinite QP, now with binary variables. 

To make the final reformulation into a MILP one must observe that 
\[ \minim{y,\bar{z}} \half y^T\bar{H}y = \minim{\bar{z}} \left(\minim{y} \half y^T\bar{H}y \right) \]
and here we can replace the inner optimisation problem with its KKT conditions in the minimisation over $\bar{z}$ and thus we arrive at
\begin{subequations} \label{eq:KKT2}
\begin{align}
  \bar{H} y + \bar{A}^T\eta + \bar{E}^T\mu = 0 \label{eq:blockedGrad}\\
  \bar{A}y - \bar{b}  - \bar{d}\bar{z} \leq 0 \\
 \bar{E}y = 0 \\
\eta^T(\bar{A}y - \bar{b} - \bar{d}\bar{z}  ) = 0 \\
  \eta \geq 0 \label{eq:KKT2-last}
\end{align}
\end{subequations}
If we now multiply \eqref{eq:blockedGrad} with $\half y^T$ from the left we have
\[ \half y^T\bar{H} y = - \half (y^T\bar{A}^T\eta + \underbrace{y^T \bar{E}^T}_{=0}\mu) = -\half \left( \bar{b} + \bar{d}\bar{z} \right)^T \eta \] 
We see that in the optimum the objective can be equivalently written as a linear term plus a bilinear term between a real variable and a binary variable. We can see that the elements of the bilinear term is either 0, if $\bar{z}^{(i)} = 0$, or euqal to $\bar{d}_i^T\eta$ when $\bar{z}^{(i)} = 1$.  Hence we can introduce a new variable, $w$, as
\[ \bar{z}^T \bar{d}^T\eta = \ones^Tw\]
where the elements $w_i$ can be modeled, yet again using the Big-M formulation, with four linear constraints
\[ -M_3(\ones-\bar{z}) \leq w - \bar{d}^T\eta \leq M_3(\ones-\bar{z}), \quad -M_3\bar{z} \leq w \leq M_3\bar{z} \]

We can now combine all the pieces together and formulate the stability problem \eqref{eq:stabilitytest} as the following MILP
\begin{subequations} \label{eq:stabilityMILP}
\begin{align}
\minimize{\eta,\mu,y,w,\bar{z},q} & - \half \bar{b}^T\eta  - \half \ones^Tw \\
\st \;\; &   \nonumber \\ 
 &  \bar{H} y + \bar{A}^T\eta + \bar{E}^T\mu = 0 \\
 & \bar{E}y = 0 \\
 & -M_3(\ones-\bar{z}) \leq w - \bar{d}^T\eta \leq M_3(\ones-\bar{z}) \\
 & -M_3\bar{z} \leq w \leq M_3\bar{z} \\
 & - M_4(1-q) \leq (\bar{A}y - \bar{b} - \bar{d}\bar{z}  ) \leq 0 \\
 & 0 \leq \eta \leq M_5q 
\end{align}
\end{subequations}

Note that we are not really interested in finding the optimum of problem \eqref{eq:stabilityMILP}, but rather find if there exist \emph{one} point, $x_k$, where the objective is less than zero. Hence we can reformulate the problem into a feasibility problem which often is much faster to solve.
\begin{subequations} \label{eq:feasibilityMILP}
\begin{align}
\text{find} \;\; & x_k \nonumber\\
\st \;\; &   \nonumber \\
& - \half \bar{b}^T\eta  - \half \ones^Tw < 0 \\
  &  \bar{H} y + \bar{A}^T\eta + \bar{E}^T\mu = 0 \\
 & \bar{E}y = 0 \\
 & -M_3(\ones-\bar{z}) \leq w - \bar{d}^T\eta \leq M_3(\ones-\bar{z}) \\
 & -M_3\bar{z} \leq w \leq M_3\bar{z} \\
 & - M_4(1-q) \leq (\bar{A}y - \bar{b} - \bar{d}\bar{z}  ) \leq 0 \label{eq:primal_q}\\
 & 0 \leq \eta \leq M_5q \label{eq:dual_q}
\end{align}
\end{subequations}

Due to the special structure of the original problem there exist a lot of structure in the problem \eqref{eq:feasibilityMILP} that should be exploited to enhance the performance of the MILP representation.

\section{Exploiting structure in the MILP}
First, let us observe that in MPC problems there are often upper and lower bounds on the variables, e.g., $u_{min} \leq u_{k+i} \leq u_{max}$. These constraints can not be fulfilled with equality at the same time and hence the two corresponding binary variables, $z_k^{(i)}$ and $z_k^{(j)}$, in \eqref{eq:compl_bin} can not be equal to one at the same time. So we can introduce the constraint
\begin{equation}
  \label{eq:binary_cond1}
  z_k^{(i)} + z_k^{(j)} \leq 1
\end{equation}
for the appropriate indices $i$ and $j$.

Let us now look at a single constraint $e_i^Tx_k + f_i^TU_k \leq b_i$ in the original MPC problem. This single constraint generates through \eqref{eq:compl_bin} the four constraints 
\begin{subequations}
\begin{align}
e_i^Tx_k + f_i^TU_k - b_i &\leq 0 \label{eq:bin_con_1a}\\
b_i - e_i^Tx_k -  f_i^TU_k &\leq m_2^{(i)}(1 - z_k^{(i)}) \label{eq:bin_con_1b}\\
- \lambda_k^{(i)} &\leq 0 \label{eq:bin_con_2a}\\
\lambda_k^{(i)} &\leq m_1^{(i)} z_k^{(i)} \label{eq:bin_con_2b}
\end{align}
\end{subequations}
in \eqref{eq:blocked_constr}.

By formulating the KKT conditions for \eqref{eq:blocked}, for each of these four constraints we have yet another binary variable, $q^{(\cdot)}$, in \eqref{eq:primal_q} and \eqref{eq:dual_q}. Since the binary variables forces a constraint to be active we can see that if $z_k^{(i)} = 0$ then \eqref{eq:bin_con_1a} and \eqref{eq:bin_con_1b} can not both be active at the same time. This means that the corresponding elements $q^{(i)}$ and $q^{(n_\lambda+i)}$ can not both be equal to one, i.e., we can constrain them as
\begin{equation}
  \label{eq:binary_cond2}
  q^{(i)} + q^{(n_\lambda+i)} \leq 1 + z_k^{(i)}
\end{equation}
Furthermore if $z_k^{(i)} = 1$ we have from \eqref{eq:bin_con_1a} and \eqref{eq:bin_con_1b} that $e_i^Tx_k+f_i^TU_k = b_i$ which gives in \eqref{eq:primal_q} and \eqref{eq:dual_q} that 
\begin{align*}
m_4^{(i)}(1-q^{(i)}) &\leq e_i^Tx_k+f_i^TU_k - b_i = 0 \\
m_4^{(n_\lambda+i)}(1-q^{(n_\lambda+i)}) &\leq b_i - e_i^Tx_k - f_i^TU_k = 0
\end{align*}
and hence we can in this case without any loss of generality constrain both $q^{(i)}$ and $q^{(n_\lambda+i)}$ to be equal to one. This is done by adding the constraint
 \begin{equation}
  \label{eq:binary_cond3}
  q^{(i)} + q^{(n_\lambda+i)} \geq 2 z_k^{(i)}
\end{equation}
The corresponding argumentation can be used for \eqref{eq:bin_con_2a} and \eqref{eq:bin_con_2b} to introduce the two additional 
constraints
\begin{align}
q^{(3n_\lambda+i)} + q^{(4n_\lambda+i)}  &\leq 1 + (1-z_k^{(i)}) \label{eq:binary_cond4} \\
q^{(3n_\lambda+i)} + q^{(4n_\lambda+i)}  &\geq 2 (1-z_k^{(i)}) \label{eq:binary_cond5}
\end{align}

The constraints \eqref{eq:binary_cond2} -\eqref{eq:binary_cond5} concerns relations in one MPC contraint, $i$, but as stated in \eqref{eq:binary_cond1} there exist relations also between different MPC constraints. This gives additional relationships between the different binary variables which we encode as the following constraints
\begin{align}
  q^{(n_\lambda+i)} &\leq 1-z_k^{(j)} \\
  q^{(2n_\lambda+i)} &\leq 1-z_k^{(j)} \\
  q^{(3n_\lambda+i)} + q^{(4n_\lambda+i)} &\geq 2z_k^{(j)} \\
  q^{(n_\lambda+j)} &\leq 1-z_k^{(i)} \\
  q^{(2n_\lambda+j)} &\leq 1-z_k^{(i)} \\
  q^{(3n_\lambda+j)} + q^{(4n_\lambda+j)} &\geq 2z_k^{(i)} \label{eq:binary_cond_x}
\end{align}

The stability analysis optimisation problem to be solved thus consist of the problem \eqref{eq:feasibilityMILP} with the additional constraints \eqref{eq:binary_cond1} and \eqref{eq:binary_cond2} - \eqref{eq:binary_cond_x}. Note that the added binary constraints are redundant in the original problem and does not effect the optimal solution. They are only added to cut of binary combinations in order to possibly increase the performance of the solver. 

Note that nowhere in the derivation of the algoritm do we use the fact that it is the MPC controller objective function that we have as candidate Lyapunov function, $V_k$. Hence we can generalize the algorithm to the use of any other positive definite $V_k$ by appropriately modifying the matrix $\bar{H}$.

\section{Examples}
In this section we will look at three examples and try to illustrate some properties and performance of the proposed algorithm.

All examples have been implemented in Matlab using YALMIP, \cite{Lofberg2004}. The MILP problems have been solved using the solver Gurobi 5.6.2, \cite{GurobiOptimization2015} and the LMI problems have been solved using MOSEK 7.1, \cite{mosek2015}.

\subsection{Sufficient but not necessary condition}
The first example is taken from \cite{Lofberg2003} where we consider the following unstable system 
\[ x_{k+1} = \begin{bmatrix} 1.216 & -0.055 \\ 0.221 & 0.9947 \end{bmatrix}x_k + \begin{bmatrix} 0.02763 \\ 0.002673 \end{bmatrix} u_k \]
which we control using an unconstrained finite time MPC controller with the objective function
\[ V_k = \minimize{u_{k+i}} \sum_{i=0}^{N-1} x_{k+i}^TQx_{k+i} + Ru_{k+i}^2 \]
and $Q = 10I$ and $R=1$. Since the controller is unconstrained can we easily calculate for which prediction horizon length, $N$, the closed loop system is stable by looking at the eigenvalues of the closed loop system matrix, $A+BL$. In the range $N = \range{1}{50}$ we calculate the eigenvalues of the closed loop system and can conclude that it is stable for $N \geq 8$. When we use the algorithm derived in section \ref{seq:the_algorithm} to test for stability we get the result that the system is stable for $N \geq 21$. Running the LMI algorithm from \cite{Primbs2001} we obtain the same results. 

Here we clearly see that the test is only a sufficient, but not necessary, condition for stability since apparently there exist a set of $8  \leq N \leq 20$ where the system is stable but the objective function does not constitute a valid Lyapunov function for the closed loop system.  

\subsection{Complexity}
To investigate the complexity of the proposed algorithm we consider a very simple example. The system is the following stable two state system
\[ x_{k+1} = \begin{bmatrix} 0.9744 & 0.0141 \\ -0.1023 & 0.9003 \end{bmatrix}x_k + \begin{bmatrix} 0.0106 \\ 0.4878 \end{bmatrix} u_k \]
This system is controlled with an input constrained MPC controller with objective function 
\[ \sum_{i=0}^{N-1} x_{k+i}^TQx_{k+i} + Ru_{k+i}^2 + x_{k+N}^TPx_{k+N} \]
with $Q = I$, $R=10$ and $P=\begin{bmatrix}10 & 1 \\1 & 2 \end{bmatrix}$ and the constraints $-5 \leq u_{k+i} \leq 5$. 

To investigate the computational complexity of the algorithm we have compared our developed MILP algorithm to the LMI algorithm developed in \cite{Primbs2001}.

Both algorithms verifies that the closed loop system is stable for all tested prediction horizon lengths but the computational time differs quite much between the two algorithms.
 
The following table shows the solvers computational time in seconds for the two methods as a function of prediction horizon.
\begin{center}
\begin{tabular}{| l || c | c | c | c | c |}
\hline
N & 2 & 4 & 6 & 8 & 10 \\ \hline \hline
LMI & 0.03 & 0.52 & 5.36 & 22.42 & 80.06 \\ \hline
MILP & 0.06 & 0.17 & 0.28 & 0.42 & 0.61\\
\hline
\end{tabular}
\end{center}
As can be seen from the table the computation time grows very rapidly with increasing prediction horizon for the LMI algorithm while it remain relatively small for the MILP approach. 

\subsection{Application to aircraft control}
Let us now consider a more realistic example taken from the aircraft industry.

\begin{figure}[hbt]
  \centering
  \includegraphics[width=0.35\textwidth]{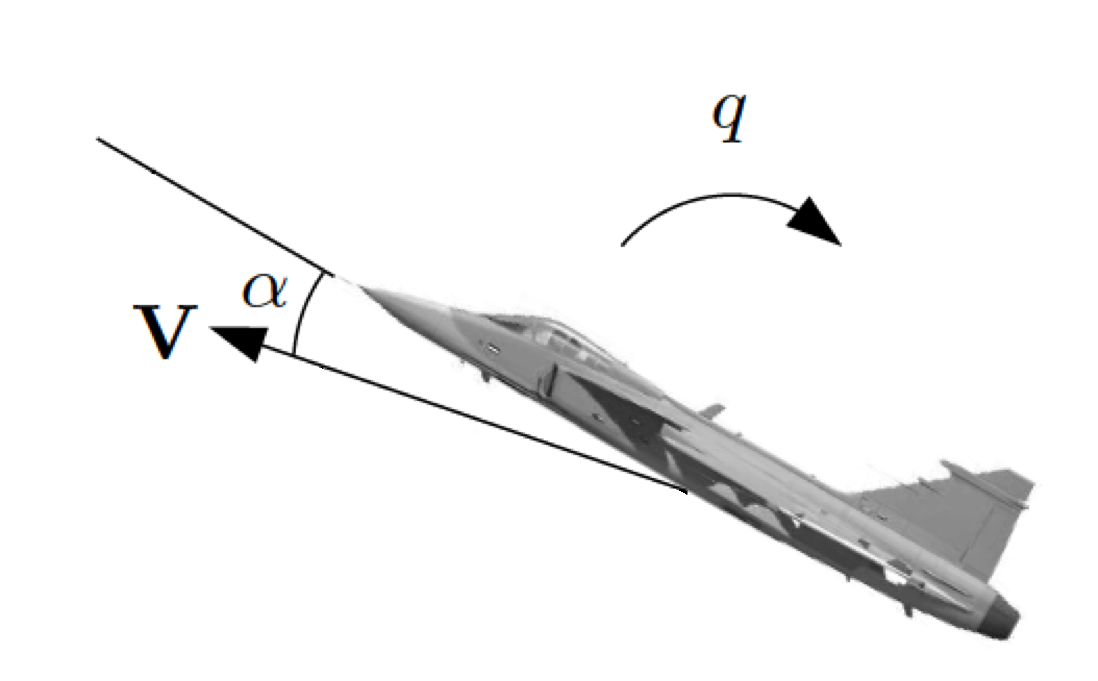}
  \caption{The two states, angle of attack ($\alpha$) and pitch rate ($q$) of the linear short period dynamic's approximation}
  \label{fig:ac}
\end{figure}

We consider the stabilisation of the so called short period dynamics, \cite{Stevens2003}, of a fighter aircraft. The short period dynamics can be approximated with a two state discrete time linear system where $x_k = [\alpha_k \quad q_k]$ is the angle of attack and pitch rate, see Figure~\ref{fig:ac}, and the input is the control surface deflection. 

The short period dynamics considered in this example are those of the ADMIRE aircraft, \cite{Forssell2005}, and we have linearised the system at Mach 0.6 and altitude 4km. This result in the following system 
\[  x_{k+1} = \begin{bmatrix} 0.9798 & 0.0158 \\ 0.1449 & 0.9787 \end{bmatrix} x_k + \begin{bmatrix} 0.0106 \\ 0.4878 \end{bmatrix} u_k \]
For this system we design a stabilising MPC controller of the form \eqref{eq:MPC} where the objective is a quadratic function
\[ \sum_{i=0}^{N-1} x_{k+i}^TQx_{k+i} + Ru_{k+i}^2 + x_{k+N}^TPx_{k+N} \]
with $Q = \begin{bmatrix} 2 & 0 \\ 0 & 0.1 \end{bmatrix}$, $R=10$ and $P$ is the associated LQ cost. The constraints are upper and lower bounds on the states and control. 
\[ \begin{bmatrix} -10 \\ -50 \end{bmatrix} \leq x_k \leq  \begin{bmatrix} 10 \\ 50 \end{bmatrix}, \quad -20 \leq u_k \leq 20 \] 
Since the system is unstable and constrained it is ''generally necessary'' to have a final state constraint, \cite{Mayne2000}. We
select the final state constraint set to the invariant set of the associated LQ controller, see e.g., \cite{Rawlings2009} for more details. 

Let us first consider the statement that it is necessary to have a terminal state constraint. The MPC controller is by construction stabilising and when we test the stability of the closed loop system with the algorithm derived in previous section we see that it indeed is stable for all $N$ in the tested range ($N = \range{2}{10}$). The question is, is it still stable if we would remove the terminal state constraint set?  

\begin{figure}[hbt]
  \centering
  \includegraphics[width=0.35\textwidth]{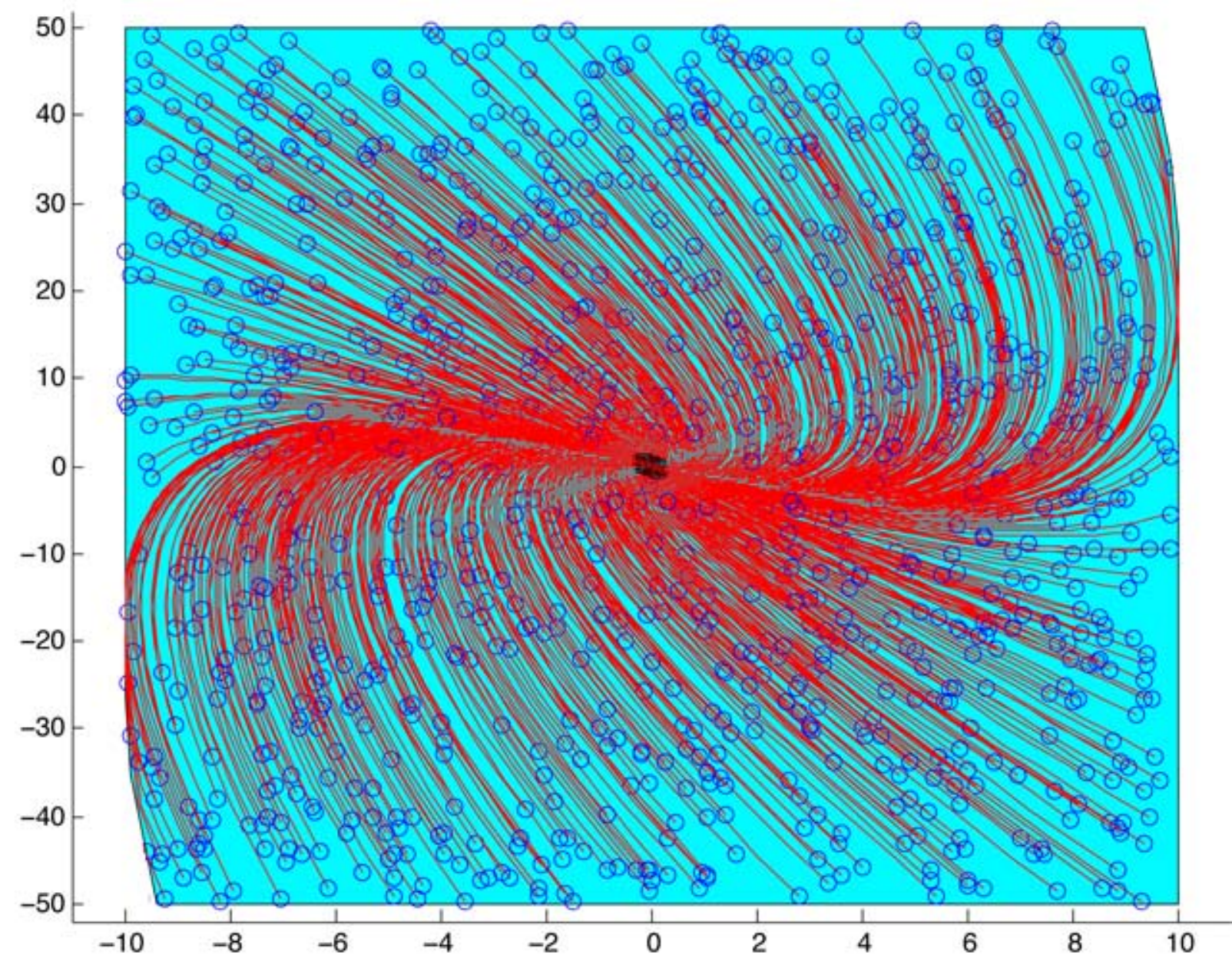}
  \caption{A simulation of the closed loop system from 1000 different initial conditions for the MPC controller without terminal constraint. The blue circles are the initial conditions, red lines are the trajectories and the black squares are the final states. The cyan polytope is the set of initially feasible states.}
  \label{fig:batch_sim}
\end{figure}

Running the algorithm again, now without the terminal constraint set in the controller, the test indicate that the controller still is stable for all tested prediction horizons. To verify this result we select 1000 random initial conditions and simulate the closed loop system. As shown in Figure \ref{fig:batch_sim} all initial conditions that are within the initially feasible set are stable.

Let us now simplify the MPC controller even more by introducing move blocking and analyse the stability of the resulting closed loop system. We select the prediction horizon $N=4$ and introduce the move blocking structure 
\[ U_k = T \hat{U}_k, \quad T = \begin{bmatrix} 1 & 0\\ 1 & 0 \\ 0 &1 \\ 0 & 1 \end{bmatrix} \] 
where $\hat{U}_k$ is the new reduced set of input signals. This blocking structure gives a controller where the sampling time of the controller output is half the internal sampling time used in the predictions.

As pointed out in \cite{Cagienard2007} the standard feasibility and stability arguments can not be used to prove stability for this move blocking strategy. Instead we aim to prove stability by using our MILP test.

\begin{figure}[hbt]
  \centering
  \includegraphics[width=0.35\textwidth]{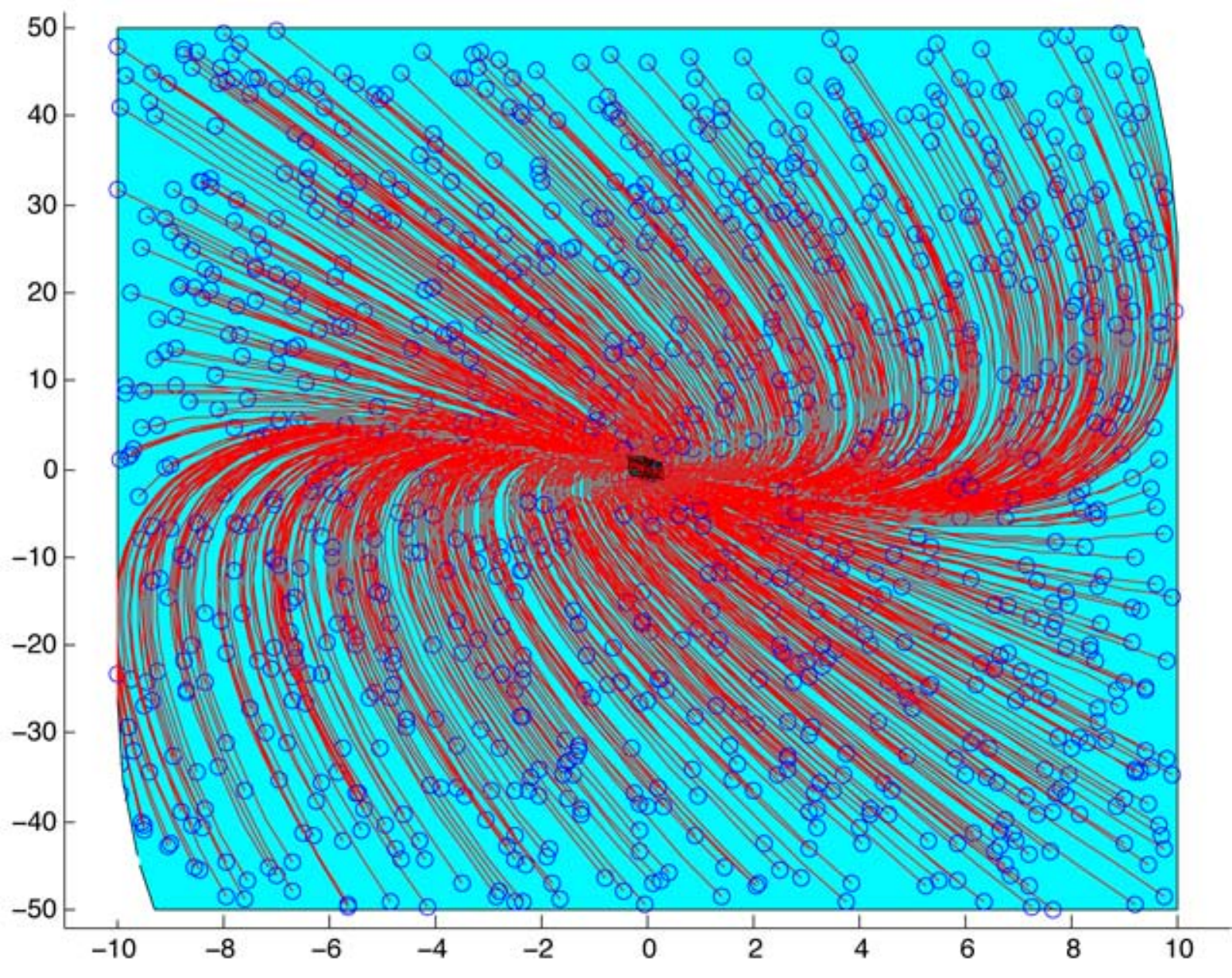}
  \caption{A simulation of the closed loop system from 1000 different initial conditions for the move blocking MPC controller.}
  \label{fig:batch_sim_mb}
\end{figure}

When applying the MILP algorithm it returns a sufficient condition certificate that the move blocking MPC design is stabilising. To verify this result we again simulate a set of 1000 random initial conditions throughout the state constraint set, see Figure~\ref{fig:batch_sim_mb}. It is clear from this figure that for all initially feasible points the move blocking MPC controller remains feasible and stabilises the system, as proven by the MILP test. 

\section{Conclusions}
In this paper we have derived an algorithm for a-postiori stability analysis of MPC controllers. 

The test allows the control engineer to verify more complex MPC designs that normally is not possible to a-priori guarantee stability for. We have shown that the stability test can be written as a Mixed Integer Linear Program which has far less computational complexity than other tests based on solving LMIs.  

As illustrated in one of the examples the test is only a sufficient test for stability and if the test fails no conclusive statements can  be made about the lack of stability of the closed loop system. 

Future work consists of further analysing the structure of the problem to reduce the computational complexity. Also extending the algorithm to a robust stability test for systems with disturbances and modelling errors. 

\vspace{6pt}

\subsubsection*{Acknowledgement}
This work has been done as a cooperation between Link\"oping University and Saab Aeronautics and is funded by the \emph{Swedish Governmental Agency for Innovation Systems} (VINNOVA) and \emph{Centrum f\"or industriell informationsteknologi} (CENIIT). 

\bibliographystyle{ieeetran}
\bibliography{library}

\end{document}